\documentclass[journal]{IEEEtran}
\usepackage{xcolor}
\usepackage{graphicx}
\usepackage{amsmath,amssymb}
\usepackage{xcolor}
\usepackage{setspace} 
\usepackage{units}
\usepackage[absolute,overlay]{textpos}
\usepackage{multirow}
\usepackage{multicol}
\usepackage{amsthm}
\usepackage{bm}
\usepackage{subfigure}
\usepackage{calligra}
\newtheorem{assum}{Assumption}
\newtheorem{definition}{Definition}

\newtheorem{theorem}{Theorem}

\title{\LARGE \bf
Event-Triggered Source Seeking Control for Nonholonomic Systems}

\author{Victor Hugo Pereira Rodrigues$^{1}$, Tiago Roux Oliveira$^{1}$ and Miroslav Krsti{\' c}$^{2}$
\thanks{This study was financed in part by the Coordena{\c c}{\~a}o de Aperfei{\c c}oamento de Pessoal de N{\'i}vel Superior – Brasil (CAPES) – Finance Code 001. The authors also acknowledge the Brazilian Funding Agencies Conselho Nacional de Desenvolvimento Cient{\'i}fico e Tecnol{\'o}gico (CNPq) and Funda{\c c}{\~a}o de Amparo {\`a} Pesquisa do Estado do Rio de Janeiro (FAPERJ).}
\thanks{$^{1}$V. H. P. Rodrigues and T. R. Oliveira are with the Department of Electronics and Telecommunication Engineering, State University of Rio de Janeiro, Rio de Janeiro -- RJ, Brazil (e-mail: victor.rodrigues@uerj.br and tiagoroux@uerj.br).}%
\thanks{$^{2}$M. Krsti{\' c} is with the Department of Mechanical and Aerospace\linebreak  Engineering at the Jacobs School of Engineering, University of California San Diego, La Jolla -- CA, USA (e-mail: mkrstic@ucsd.edu).}%
}

\begin{document}

\maketitle
\thispagestyle{empty}
\pagestyle{empty}

\begin{abstract}
This paper introduces an  event-triggered source seeking control (ET-SSC) for autonomous vehicles modeled as the nonholonomic unicycle. The classical source seeking control is enhanced with static-triggering conditions to enable aperiodic and less frequent updates of the system's input signals, offering a resource-aware control design. Our convergence analysis is based on time-scaling combined with Lyapunov and averaging theories for systems with discontinuous right-hand sides. ET-SSC ensures exponentially stable behavior for the resulting average system, leading to practical asymptotic convergence to a small neighborhood of the source point. We guarantee the avoidance of Zeno behavior by establishing a minimum dwell time to prevent infinitely fast switching. The performance optimization is aligned with classical continuous-time source seeking algorithms while balancing system performance with actuation resource consumption. Our ET-SSC algorithm, the first of its kind, allows for arbitrarily large inter-sampling times, overcoming the limitations of classical sampled-data implementations for source seeking control. 
\end{abstract}

\section{INTRODUCTION}

In the modern era of network science, researchers are increasingly focused on reducing costs by developing fast and reliable communication strategies, particularly in scenarios where the plant and controller are not physically connected or are located in different geographical regions. Networked control systems offer significant advantages in terms of installation and maintenance costs \cite{ZHGDDYP:2020}. However, a major drawback of these systems is the high traffic congestion they generate, which can result in transmission delays and packet dropouts—data losses that occur while information is transmitted through the network \cite{HNX:2007}. These challenges are closely linked to limited resources and the constrained bandwidth of available communication channels. To address this issue, event-triggered controllers can be employed as an effective solution.

Event-Triggered (ET) control performs the control task on a non-periodic basis, responding to a triggering condition that is defined as a function of the plant's state. In addition to ensuring asymptotic stability \cite{T:2007}, this approach minimizes control effort by updating the control signal and transmitting data only when the error between the current state and the equilibrium set surpasses a threshold that could lead to instability \cite{BH:2013}. Early contributions to resource-aware control design include the development of digital computer designs \cite{s9b}, event-based PID controllers \cite{s9}, and event-based controllers for stochastic systems \cite{s8}. Many studies have extended event-based control to networked systems with high levels of complexity, addressing both linear \cite{s1,HJT:2012,s5} and nonlinear systems \cite{T:2007,APDN:2016}.

Many engineering challenges encountered in industries such as network virtualization, software-defined networks, cloud computing, the Internet of Things, context-aware networks, green communications, and security can be modeled using a event-triggered approach. These problems often involve distributed connectivity through networks, where resources are shared and utilized effectively \cite{Basar:2019, AB:2011}. As a result, the need for real-time optimization to enhance these engineering processes is both evident and crucial. While there are existing studies addressing these topics \cite{Mazo_Tabuada, Wang_Lemmon, Johansson}, the application of extremum-seeking feedback in this context has yet to be explored \cite{KW:2000}.

Despite the extensive body of literature on ET, as well as the established results of extremum seeking (ES) for static and general nonlinear dynamic systems in continuous time \cite{KW:2000}, the theoretical developments of ES have extended well beyond these systems. These advancements now encompass discrete-time systems \cite{CKAL:2002}, stochastic systems \cite{MK:2009, LK:2010}, multivariable systems \cite{GKN:2012}, noncooperative games \cite{FKB:2012}, time delays \cite{OKT:2017, ZFO:2023}, and even broader classes of infinite-dimensional systems governed by partial differential equations (PDEs) \cite{OK:2022_book_golden}. However, until recently, no work directly addressed the integration of ES with ET---this gap was only filled by \cite{AUT:2025}, which developed multi-variable ES algorithms based on perturbation (averaging) estimates of the model through ET control. Although progress has been made, the challenge of simultaneously addressing \textit{source seeking control} \cite{SCL_2007,source1,source2,source3,source4,source5,source6} in GPS-denied environments and integrating event-triggered versions of extremum seeking control remains an open problem.

Source seeking control (SSC) addresses the problem of locating the source of a scalar signal using an autonomous vehicle modeled as a non-holonomic unicycle, which is equipped with a sensor to detect the scalar signal but lacks the ability to sense either the position of the source or its own location. The authors in \cite{SCL_2007} assume that the signal field is strongest at the source and decays as it moves away from it, though the vehicle does not have access to the functional form of the field. To guide the vehicle toward the source, they employ extremum seeking to estimate the gradient of the field in real time. The goal is to steer the vehicle towards the point where the gradient is zero, which corresponds to the location of the source (the maximum of the field). The vehicle follows a periodic forward-backward motion, implementable on mobile robots and some underwater vehicles, but not on aircraft. The forward velocity includes a tunable bias term, which, when combined with extremum seeking, enables the vehicle to drift towards the source. 

In this paper, we extend previous designs from the multi-input-single-output ES scenario to SSC framework, through a ET-SSC approach, enabling arbitrarily large inter-event times and overcoming the maximum sampling constraints of sampled-data control. We design static triggering conditions that guarantee the absence of the \emph{Zeno phenomenon} by establishing a minimal dwell time between successive event triggers. This ensures Zeno-free behavior in both the average and original systems. Crucially, we prove that ET-SSC achieves local exponential stability in the average sense and ensures local asymptotic convergence of the unknown parameters to a neighborhood of the source. Our approach employs continuous dither signals and leverages averaging theory for systems with discontinuous right-hand sides \cite{P:1979}. The non-periodic nature of ETC does not hinder the application of Plotnikov’s averaging results \cite{P:1979} since the perturbation-probing signals in SSC remain periodic. The stability analysis of ET-SSC relies on time-scaling, averaging theory for systems with discontinuous right-hand sides, and Lyapunov’s direct method. In summary, our work presents the first \emph{systematic} framework for event-triggered source seeking based on classical periodic-perturbation techniques \cite{KW:2000}.

\section{PROBLEM FORMULATION}

We consider a unicycle-type mobile robot equipped with a sensor, which may be placed directly at the vehicle's center or at some distance \( r \) from it \cite{SCL_2007}. 
Based on this set-up, the motion of the robot center follows the equations:
\begin{align}
    \frac{dx}{dt}(t) &= v(t) \cos (\theta(t)), \label{eq:dotX_ETSSC_v1} \\
    \frac{dy}{dt}(t) &= v(t) \sin (\theta(t)), \label{eq:dotY_ETSSC_v1} \\
    \frac{d\theta}{dt}(t) &= \omega(t), \label{eq:dotTheta_ETSSC_v1}
\end{align}
where $x(t), y(t) \in \mathbb{R}$ represents the central coordinates of the vehicle, $\theta(t)\in \mathbb{R}$ is its heading angle, and $v(t), \omega(t) \in \mathbb{R}$ are the control inputs for the linear and angular speeds, respectively. 

The aim of the vehicle is to locate a source that emits a signal that decreases as the distance from the source increases. We assume that this signal field follows an unknown nonlinear function 
\begin{align}
    Q(x, y,\theta)&\mathbb{=}Q^{\ast}\mathbb{-}\frac{1}{2}(x(t)\mathbb{-}x^{\ast})^2\mathbb{-}\frac{1}{2}(y(t)\mathbb{-}y^{\ast})^2\mathbb{-}\frac{1}{2}(\theta(t)\mathbb{-}\theta^{\ast})^2\!\!,  \label{eq:Q_ETSSC_v1}  
\end{align}
which has a single local maximum $Q^* = Q(x^{\ast},y^{\ast},\theta^{\ast})$ at point $(x^{\ast},y^{\ast},\theta^{\ast})$, representing the source's location. Our control objective is to design a control strategy that guides the vehicle to this maximum using only the signal values measured at the sensor, without prior knowledge of the function's parameters, according to the following assumption.

\begin{assum}\label{assumption1}
The unique optimizer vector $[x^{\ast}\,, y^{\ast}\,,\theta^{\ast}]^{\top} \in \mathbb{R}^{3}$ and the scalar $Q^{\ast} \in \mathbb{R}$ are unknown parameters of the nonlinear map (\ref{eq:Q_ETSSC_v1}).
\end{assum}

In this paper, we particularly propose an Event-Triggered Source Seeking Control (ET-SSC) that executes control actions only when a desired condition is specified, based on the estimated gradient. This event-triggering mechanism allows the system to operate aperiodically, reducing computational load and resource usage by avoiding unnecessary control updates. The goal of this strategy is to drive the robot pose $[x(t), y(t), \theta(t)]$ to converge to the source location $[x^{\ast}, y^{\ast}, \theta^{\ast}]$. 

\subsection{Linear Velocity $v(t)$ and the Angular Velocity $\omega(t)$}

By assuming the sensor is placed at the robot’s center, both the linear velocity $v(t)$ and the angular velocity $\omega(t)$ are adjusted according to the following tuning laws: 
\begin{align}
     v(t)&=\cos(\theta(t))[a_{1}\omega_{1}\cos(\omega_{1}t)+u_{1}(t)]+ \nonumber \\
     &\quad+\sin(\theta(t))[a_{2}\omega_{2}\sin(\omega_{2}t)+u_{1}(t)]\,, \label{eq:vt_ETSSC_v1} \\
     \omega(t)&=\frac{a_{3}\omega_{3}}{2}\cos(\omega_{3}t)+u_{2}(t) \label{eq:wt_ETSSC_v1}\,,
\end{align}
where $u_{1}(t), u_{2}(t) \in \mathbb{R}$ are ET-SSC laws to be designed,
$a_{i} >0$, for $i \in \{1,2,3\}$, are small dither's amplitudes and $\omega_{i}$'s are the probing frequencies. 
%
We make the following
assumption about them. 
\begin{assum}\label{assum:w}
The probing frequencies satisfy
\begin{align}
\omega_1 = \omega_2 = 2\omega_3\,, \label{eq:w1_w2}
\end{align}
where $\omega_{3}$ is a positive constant.
\end{assum} 

By substituting equation (\ref{eq:vt_ETSSC_v1}) into equations (\ref{eq:dotX_ETSSC_v1}) and (\ref{eq:dotY_ETSSC_v1}), and replacing equation~(\ref{eq:wt_ETSSC_v1}) into equation (\ref{eq:dotTheta_ETSSC_v1}), the system dynamics described by equations (\ref{eq:dotX_ETSSC_v1})--(\ref{eq:dotTheta_ETSSC_v1}) can be written as
\begin{align}
    &\frac{dx}{dt}(t) = \frac{a_{1}\omega_{1}}{2}\cos(\omega_{1}t)\mathbb{+}\frac{1}{2}(1+\sin (2\theta(t))\mathbb{+}\cos (2\theta(t)))u_{1}(t) \nonumber \\
    &\mathbb{+}\left(\frac{a_{1}\omega_{1}}{2}\cos(\omega_{1}t)\mathbb{+}\frac{a_{2}\omega_{2}}{2}\sin(\omega_{2}t)\right)(\sin (2\theta(t))\mathbb{+}\cos (2\theta(t)))  \label{eq:dotX_ETSSC_v2} \\
   &\frac{dy}{dt}(t) = \frac{a_{2}\omega_{2}}{2}\sin(\omega_{2}t)\mathbb{+}\frac{1}{2}(1+\sin (2\theta(t))\mathbb{-}\cos (2\theta(t)))u_{1}(t) \nonumber \\
    &\mathbb{+}\left(\frac{a_{1}\omega_{1}}{2}\cos(\omega_{1}t)\mathbb{+}\frac{a_{2}\omega_{2}}{2}\sin(\omega_{2}t)\right)(\sin (2\theta(t))\mathbb{-}\cos (2\theta(t))) \,, \label{eq:dotY_ETSSC_v2} \\
    &\frac{d\theta}{dt}(t) = \frac{a_{3}\omega_{3}}{2}\cos(\omega_{3}t)+u_{2}(t)\,. \label{eq:dotTheta_ETSSC_v2}
\end{align}

Hence, by defining, 
\begin{align}
    &\frac{d\hat{x}}{dt}(t) = \frac{1}{2}(1+\sin (2\theta(t))\mathbb{+}\cos (2\theta(t)))u_{1}(t) \nonumber \\
    &\mathbb{+}\left(\frac{a_{1}\omega_{1}}{2}\cos(\omega_{1}t)\mathbb{+}\frac{a_{2}\omega_{2}}{2}\cos(\omega_{2}t)\right)(\sin (2\theta(t))\mathbb{+}\cos (2\theta(t))), \label{eq:dotHatX_ETSSC_v1} \\
   &\frac{d\hat{y}}{dt}(t) = \frac{1}{2}(1+\sin (2\theta(t))\mathbb{-}\cos (2\theta(t)))u_{1}(t) \nonumber \\
    &\mathbb{+}\left(\frac{a_{1}\omega_{1}}{2}\cos(\omega_{1}t)\mathbb{+}\frac{a_{2}\omega_{2}}{2}\cos(\omega_{2}t)\right)(\sin (2\theta(t))\mathbb{-}\cos (2\theta(t))), \label{eq:dotHatY_ETSSC_v1} \\
    &\frac{d\hat{\theta}}{dt}(t) = u_{2}(t)\,, \label{eq:dotHatTheta_ETSSC_v1}
\end{align}
it is possible to write
\begin{align}
    x(t)&=\hat{x}(t)+\frac{a_{1}}{2}\sin(\omega_{1}t)\,, \label{eq:x_ETSSC_v1} \\
    y(t)&=\hat{y}(t)-\frac{a_{2}}{2}\cos(\omega_{2}t)\,, \label{eq:y_ETSSC_v1} \\
    \theta(t)&=\hat{\theta}(t)+\frac{a_{3}}{2}\sin(\omega_{3}t)\,. \label{eq:theta_ETSSC_v1}
\end{align}

\subsection{Estimation Error Dynamics}

We define the estimation errors as 
\begin{align}
    \tilde{x}(t)&=\hat{x}(t)-x^{\ast}\,, \label{eq:tildeX_ETSSC_v1} \\
    \tilde{y}(t)&=\hat{y}(t)-y^{\ast}\,, \label{eq:tildeY_ETSSC_v1} \\
    \tilde{\theta}(t)&=\hat{\theta}(t)-\theta^{\ast}\,. \label{eq:tildeTheta_ETSSC_v1}
\end{align}
Now, by using (\ref{eq:tildeX_ETSSC_v1})--(\ref{eq:tildeTheta_ETSSC_v1}),  it is possible to rewrite (\ref{eq:x_ETSSC_v1})--(\ref{eq:theta_ETSSC_v1}) as
%
\begin{align}
    x(t)&=\tilde{x}(t)+x^{\ast}+\frac{a_{1}}{2}\sin(\omega_{1}t)\,, \label{eq:x_ETSSC_v2} \\
    y(t)&=\tilde{y}(t)+y^{\ast}-\frac{a_{2}}{2}\cos(\omega_{2}t)\,, \label{eq:y_ETSSC_v2} \\
    \theta(t)&=\tilde{\theta}(t)+\theta^{\ast}+\frac{a_{3}}{2}\sin(\omega_{3}t)\,. \label{eq:theta_ETSSC_v2}
\end{align}
Then, by defining the state vector, the source vector, and estimation error vector as  
\begin{align}
    q(t) &:= [x(t),\ y(t),\ \theta(t)]^\top \,,  \label{eq:q_ETSSC_v1} \\
    q^\ast &:= [x^\ast,\ y^\ast,\ \theta^\ast]^\top  \,,  \label{eq:qAst_ETSSC_v1} \\
    \tilde{q}(t) &:=  [\tilde{x}(t)\,, \tilde{y}(t)\,, \tilde{\theta}(t)]^\top \,, \label{eq:tildeQ_ETSSC_v1}
\end{align} 
and the dither vector  
\begin{align}
S(t) &:= \left[\frac{a_1}{2}\sin(\omega_1 t),\ -\frac{a_2}{2}\cos(\omega_2 t),\ \frac{a_3}{2}\sin(\omega_3 t)\right]^\top, \label{eq:S_v1}
\end{align}
we can rewrite (\ref{eq:x_ETSSC_v2})--(\ref{eq:theta_ETSSC_v2}) in the vector form
\begin{align}
q(t) - q^\ast = \tilde{q}(t) + S(t)\,. \label{eq:q_ETSSC_v2}
\end{align}  

By the assumption of a small deviation of orientation angle, {\it i.e.}, $2\tilde{\theta}(t)\ll 1$, one has $\sin(2\tilde{\theta}(t)) \approx 2\tilde{\theta}(t)$ and $\cos(2\tilde{\theta}(t)) \approx 1$ and neglecting quadratic terms $\theta(t)u_{1}(t)$ as well as  $\omega_1 = \omega_2 = 2\omega_3$, the resulting dynamics governing (\ref{eq:tildeQ_ETSSC_v1}) is simplified to the linear dynamics given by 
\begin{align}
\frac{d\tilde{q}}{dt}(t)=A(t)\tilde{q}(t)+B(t)u(t)+\Delta(t)\,, \label{eq:dotq_ETSSC_v1}
\end{align}
with the input signal $u(t):= [u_{1}(t)\,,u_{2}(t)]^\top$, and the time-varying matrices $A(t), B(t)$ as well as the equivalent disturbance $\Delta(t)$ such that  
\begin{align}
A(t)&:=A+\Delta A(t)\,, \label{eq:A_ETSSC_v1} \\ 
A&=\begin{bmatrix} 
        0 & 0 & \frac{\sqrt{2}}{2}a_{1}\omega_{3}\cos\left(2\theta^{\ast}+\frac{\pi}{4}\right)J_{2}(a_{3}) \\
        0 & 0 & \frac{\sqrt{2}}{2}a_{1}\omega_{3}\cos\left(2\theta^{\ast}-\frac{\pi}{4}\right)J_{2}(a_{3}) \\
        0 & 0 & 0
      \end{bmatrix}\,, \label{eq:A_ETSSC_v2} \\  
B(t)&:= B+\Delta B(t) \,, \label{eq:B_ETSSC_v1} \\
B&=\begin{bmatrix} 
        \frac{1}{2}+\frac{\sqrt{2}}{2}\cos\left(2\theta^{\ast}-\frac{\pi}{4}\right)J_{0}(a_{3}) & 0 \\
        \frac{1}{2}-\frac{\sqrt{2}}{2}\cos\left(2\theta^{\ast}+\frac{\pi}{4}\right)J_{0}(a_{3}) & 0 \\
        0 & 1
      \end{bmatrix}     \,. \label{eq:B_ETSSC_v2}\\
\Delta(t)&:= \bar{\Delta}+\bar{\Delta}(t) \,,  \label{eq:Delta_ETSSC_v1} \\
\bar{\Delta}&\mathbb{=} \begin{bmatrix}\frac{\sqrt{2}}{2}a_{1}\omega_{3}\cos\left(2\theta^{\ast}-\frac{\pi}{4}\right)J_{2}(a_{3}) \\ -\frac{\sqrt{2}}{2}a_{1}\omega_{3}\cos\left(2\theta^{\ast}-\frac{\pi}{4}\right)J_{2}(a_{3})\\0 \end{bmatrix}\,, \label{eq:barDelta}
\end{align}
where $\Delta A(t) \in \mathbb{R}^{3 \times 3}$, $\Delta B(t) \in \mathbb{R}^{3 \times 2}$ and $\bar{\Delta} (t) \in \mathbb{R}^{3}$ have zero mean over time, {\it i.e.}, they do not introduce any bias.

In this specific case, the matrices $A(t)$ and $B(t)$ contain time-dependent components involving expressions such as $\cos(a_3\sin(\omega_3 t))$ and $\sin(a_3\sin(\omega_3 t))$, which can be expanded into series of Bessel functions of the first kind \cite{AS:1964}:
\begin{small}
\begin{align}
    \cos(a_3 \sin(\omega_3 t)) &= J_0(a_3) + 2 \sum_{m=1}^{\infty} J_{2m}(a_3) \cos(2m \omega_3 t), \label{eq:besse1} \\
    \sin(a_3 \sin(\omega_3 t)) &= 2 \sum_{m=0}^{\infty} J_{2m+1}(a_3) \sin((2m+1) \omega_3 t), \label{eq:besse2}\\
    J_{m}(a_{3}) &= \frac{1}{\pi}\int_{0}^{\pi}\cos(a_{3}\sin(\omega_{3}t)-m\omega_{3}t)d(\omega_{3}t)\,. \label{eq:besse3}
\end{align}
\end{small}

\subsection{Gradient Estimate and its Dynamics}

By plugging (\ref{eq:q_ETSSC_v1}) into (\ref{eq:Q_ETSSC_v1}), and  neglecting quadratic terms in $\tilde{q}(t)$ for a local analysis, the map $Q$ can rewritten with respect to the estimation error as 
%
%
\begin{align}
    Q(t,\tilde{q})&=Q^{\ast}-\frac{a_{1}^{2}}{16}-\frac{a_{2}^{2}}{16}-\frac{a_{3}^{2}}{16}+\frac{a_{1}^{2}}{16}\cos(2\omega_{1}t)+ \nonumber \\
    &\quad-\frac{a_{2}^{2}}{16}\cos(2\omega_{2}t)+\frac{a_{3}^{2}}{16}\cos(2\omega_{3}t)+ \nonumber \\
    &\quad-\frac{a_{1}}{2}\sin(\omega_{1}t)\tilde{q}_{1}(t)+\frac{a_{2}}{2}\cos(\omega_{2}t)\tilde{q}_{2}(t)+ \nonumber \\
    &\quad-\frac{a_{3}}{2}\sin(\omega_{3}t)\tilde{q}_{3}(t) \,.\label{eq:Q_ETSSC_v3}  
\end{align}
The gradient estimate is given by 
\begin{align}
\hat{G}(t)=M(t)Q(t,\tilde{q})\,, \label{eq:hatG_ETSSC_v1}
\end{align}
with the demodulation vector
\begin{align}
M(t)&=\left[-\frac{4}{a_{1}}\sin\left(\omega_1 t\right)\,,\frac{4}{a_{2}}\cos\left(\omega_2 t\right)\,,-\frac{4}{a_{3}}\sin\left(\omega_3 t\right)\right]^{\top} \,. \label{eq:M_ETSSC_v1}
\end{align} 
Thus, plugging (\ref{eq:Q_ETSSC_v3}) and (\ref{eq:M_ETSSC_v1}) into (\ref{eq:hatG_ETSSC_v1}), the gradient estimate is given by
\begin{align}
\hat{G}(t)&=\mbox{\calligra H}~~(t)\tilde{q}(t)+\hat{\Delta}(t)\,, \label{eq:hatG_ETSSC_v2}
\end{align}
where
\begin{align}
\mbox{\calligra H}~~(t) &:= I_{3}+\Delta\mbox{\calligra H}~~(t)\,, \label{eq:calligraH} 
\end{align}
with a time-varying disturbance $\hat{\Delta}(t) \in \mathbb{R}^{3}$, as well as a time-varying matrix $\Delta\mbox{\calligra H}~~(t)$ characterized by having zero mean over time, and $I_{3} \in \mathbb{R}^{3\times 3}$ denoting a identity matrix.

By using (\ref{eq:dotq_ETSSC_v1}), the time-derivative of (\ref{eq:hatG_ETSSC_v2}) is given by 
\begin{align}
\frac{d\hat{G}(t)}{dt}&=\mathcal{A}(t)\tilde{q}(t) +\mathcal{B}(t)u(t)+\delta(t)\,, \label{eq:dhatG_ETSSC_v_1} \\
\mathcal{A}(t)&:=\mbox{\calligra H}~~(t)A(t)+\frac{d \mbox{\calligra H}}{dt}(t)\,, \\
\mathcal{B}(t)&:=\mbox{\calligra H}~~(t)B(t) \,, \\
\delta(t)&:=\mbox{\calligra H}~~(t)\Delta(t)+\frac{d \hat{\Delta}}{dt}(t)\,.
\end{align}

\subsection{State Feedback based on Gradient Estimate}

For all $t\geq 0$, the  continuous-time feedback law
\begin{align}
u(t)=-K\hat{G}(t) \,, \quad \forall t\geq 0 \label{eq:U_continuous}
\end{align}
is a stabilizing controller for the average version of (\ref{eq:dhatG_ETSSC_v_1}), since the gain $K$ is chosen such that $A-BK$ is Hurwitz. 
%

Our goal is to design a stabilizing controller for the closed-loop system (\ref{eq:dhatG_ETSSC_v_1}) in a sampled-and-hold fashion \cite{APDNH:2018} by  emulating the continuous-time control law (\ref{eq:U_continuous}). Here, the control law is only updated for a given  sequence of time instants $(t_{k})_{k\in\mathbb{N}}$ defined by an event-generator that preserves stability and robustness. More precisely, the execution of the control task is orchestrated by a monitoring mechanism that invokes control updates when the difference between the current value of the gradient estimate and its previously computed value  at time $t_k$ becomes too large \cite{AUT:2025,HJT:2012}.

The following assumption is additionally considered throughout the paper.

\begin{assum}\label{assumption2}
The matrix product $A-BK$ is Hurwitz such that for any given $Q=Q^{\top}>0$ there exists a $P=P^{T}>0$ that satisfies the Lyapunov equation
\begin{align}
(A-BK)^{\top}P+P(A-BK)=-Q \,. \label{eq:LyapEq}
\end{align}
%
Although the matrices $A$, $B$, $P$ and $Q$ are unknown, it is possible to establish the existence of a known positive constant \( \alpha \) such that
\begin{align}
\alpha > \frac{2\|P(A - BK)\|}{\lambda_{\min}(Q)} \,, \label{eq:alpha}
\end{align}
based solely on the design specifications. This inequality holds under the assumption that the closed-loop system matrix $A-BK$ is Hurwitz and that a Lyapunov function exists satisfying the standard Lyapunov inequality.

\end{assum}

\subsection{Emulation of the Continuous-Time Extremum Seeking}

Defining the control input for all $t\in[t_{k},t_{k+1})$,  $k\in\mathbb{N}$, as
\begin{align}
u_k=-K\hat{G}(t_{k}) \,, \label{eq:U_event}
\end{align}
we introduce the error vector, which is related to the deviation of the gradient estimate as 
\begin{align}
e(t):=\hat{G}(t_{k})-\hat{G}(t) \,, \quad \forall t \in \lbrack t_{k}\,, t_{k+1}) \,, \quad k\in \mathbb{N} \,. \label{eq:e_event}
\end{align}
Thus, using  (\ref{eq:e_event}) and  (\ref{eq:U_event}), for all $t\in[t_{k}\,, t_{k+1})$,  the $i$-th event-triggered control law is rewritten as  
\begin{align}
u(t)=-K\hat{G}(t)-Ke(t) \,, \quad \forall t \in \lbrack t_{k}\,, t_{k+1}) \,, \quad k\in \mathbb{N} \,.  \label{eq:U_event2}
\end{align}

Now, using the gradient estimate (\ref{eq:hatG_ETSSC_v2}) and the event-triggered control law (\ref{eq:U_event2}), one arrives at the 
closed-loop representation of 
%
(\ref{eq:dotq_ETSSC_v1}) and (\ref{eq:dhatG_ETSSC_v_1}) with respect to the error vector  (\ref{eq:e_event}) and the time-varying disturbances
\begin{align}
\frac{d\tilde{q}}{dt}(t)&=\left(A(t)\!-\!B(t)K\mbox{\calligra H}~~(t)\right)\tilde{q}(t)\!-\!B(t)Ke(t)\!+\!\tilde{\delta}(t)\,, \label{eq:dotq_ETSSC_v2} \\
\tilde{\delta}(t)&:=-B(t)K\hat{\Delta}(t)+\Delta(t)\,, \label{eq:tildeDelta_ETSSC_v1}\\
\frac{d\hat{G}}{dt}(t)&=\mathcal{A}(t)\tilde{q}(t) -\mathcal{B}(t)K\hat{G}(t)-\mathcal{B}(t)Ke(t)+\delta(t)\,, \label{eq:dhatG_ETSSC_v_2} 
\end{align}
for all  $ t \in \lbrack t_{k}\,, t_{k+1}) \,,  ~k\in \mathbb{N}$.

The closed-loop system described by (\ref{eq:dotq_ETSSC_v2}) and (\ref{eq:dhatG_ETSSC_v_2}) highlights a crucial point: while the product $A(t)-B(t)K\mbox{\calligra H}~~(t)$ on averaging sense results in a Hurwitz matrix, the convergence to the equilibrium $\tilde{q}\equiv0$ and $\hat{G}\equiv0$ is not guaranteed due to the presence of the error vector $e(t)$ and the time-varying terms $\tilde{\delta}(t)$ and $\delta(t)$. However, the such a system is Input-to-State Stable (ISS), in the average sense, concerning the error vector $e(t)$ and such time-varying disturbances.

In the next section, we introduce a static event-triggering mechanism for source seeking control (SSC), as outlined in Definitions \ref{def:staticEvent}. This mechanism represents a fusion of event-triggered
data transmission with a source-seeking control system \cite{AUT:2025}.

\subsection{Event-Triggered Control Design}

The next definition 
employs 
\begin{align}
\Xi(\hat{G},e)&=\sigma\|\hat{G}(t)\|-\alpha(\|e(t)\|\mathbb{+}a_{1}\omega_{3}\left|J_{2}(a_{3})\right|)	\,, \label{eq:Xi_event_1}
\end{align}
with $\sigma \in (0,1)$ being a  parameter of the  static event-trigger to be designed. The  mapping $\Xi(\hat{G},e)$ is defined to properly re-compute the control law (\ref{eq:U_event}) and update the ZOH actuator such that the asymptotic stability of the closed-loop system is achieved \cite{HJT:2012}.

\begin{definition}[\small{Static Triggering Condition}] \label{def:staticEvent}
Let $\Xi(\hat{G},e)$ in (\ref{eq:Xi_event_1}) be the nonlinear mapping  and $K$  the control gain in (\ref{eq:U_event}). The event-triggered controller with static-triggering condition consists of two components:
\begin{enumerate}
	\item A set of increasing sequence of time $I=\{t_{0}\,, t_{1}\,, t_{2}\,,\ldots\}$ with $t_{0}=0$ generated under the following rules:
		\begin{itemize}
 			\item If $\left\{t \in\mathbb{R}^{+}: t>t_{k} ~ \wedge ~ \Xi (\hat{G},e) < 0\right\} = \emptyset$, then the set of the times of the events is $I=\{t_{0}\,, t_{1}\,, \ldots, t_{k}\}$.
			\item If $\left\{t \in\mathbb{R}^{+}: t>t_{k} ~ \wedge ~ \Xi (\hat{G},e) < 0\right\} \neq \emptyset$, the next event time is given by
				\begin{align}
					t_{k+1}&=\inf\left\{t \in\mathbb{R}^{+}: t>t_{k} ~ \wedge ~ \Xi (\hat{G},e) <0 \right\}\,, \label{eq:tk+1_event}
				\end{align}
				\!\!\!consisting of the static event-triggering mechanism.
		\end{itemize}
	\item A feedback control action  \eqref{eq:U_event} updated at time $t_k,$ for all $t\in[t_{k},t_{k+1})$, $k\in\mathbb{N}$.
\end{enumerate}  
\end{definition}

Now, we  introduce a suitable time scale to carry out  the stability analysis of the closed-loop system. 

\subsection{Rescaling of Time}

From (\ref{eq:w1_w2}), one can notice that the dither frequencies in (\ref{eq:vt_ETSSC_v1}), (\ref{eq:wt_ETSSC_v1})  and (\ref{eq:M_ETSSC_v1}), as well as their combinations have $\omega_{3}$ as fundamental frequency. Therefore, there exists a time period $T$ such that
\begin{align}
T&= \frac{2\pi}{\omega_{3}}\,, \label{eq:T}
\end{align}
such that it is possible to define the time scale for the dynamics (\ref{eq:dotq_ETSSC_v2}) and (\ref{eq:dhatG_ETSSC_v_2}) with the transformation 
\begin{align} \label{namadrugacomaltashoras}
\bar{t}=\omega_{3} t\,.
\end{align}

Hence, the system (\ref{eq:dotq_ETSSC_v2}) and (\ref{eq:dhatG_ETSSC_v_2}) can be rewritten as, $\forall t \in \lbrack t_{k}\,, t_{k+1}) \,,  ~k\in \mathbb{N}$: 
\begin{align}
\frac{d\tilde{q}}{d\bar{t}}(\bar{t})&=\frac{1}{\omega_{3}}\left(A(\bar{t})-B(\bar{t})K\mbox{\calligra H}~~(\bar{t})\right)\tilde{q}(\bar{t})+ \nonumber \\
&\quad-\frac{1}{\omega_{3}}B(\bar{t})Ke(\bar{t})+\frac{1}{\omega_{3}}\tilde{\delta}(\bar{t})\,, \label{eq:dotq_ETSSC_v3} \\
\frac{d\hat{G}}{d\bar{t}}(\bar{t})&=\frac{1}{\omega_{3}}\mathcal{A}(\bar{t})\tilde{q}(\bar{t}) -\frac{1}{\omega_{3}}\mathcal{B}(\bar{t})K\hat{G}(\bar{t}) \nonumber \\
&\quad-\frac{1}{\omega_{3}}\mathcal{B}(\bar{t})Ke(\bar{t})+\frac{1}{\omega_{3}}\delta(\bar{t})\,. \label{eq:dhatG_ETSSC_v_3} 
\end{align}

In the next section, the analysis of the closed-loop system will begin using the averaging method. This procedure aims to obtain an approximate description of the system's dynamic behavior over time by eliminating the fast components associated with the fundamental frequency and focusing on the slow dynamics.

\subsection{Closed-Loop Average System}

From (\ref{eq:dotq_ETSSC_v3}) and (\ref{eq:dhatG_ETSSC_v_3}), the corresponding average system can be obtained using \cite{P:1979}. Note that the right-hand side of (\ref{eq:dotq_ETSSC_v3}) and (\ref{eq:dhatG_ETSSC_v_3}) are aperiodically discontinuous in the state dynamics since the triggering events are not periodic. However,  the system is still periodic in time, which justifies the averaging theory \cite{P:1979}. 

Defining the augmented state as follows
\begin{align}
X(\bar{t}):=\begin{bmatrix} \tilde{q}(\bar{t}) \\ \hat{G}(\bar{t})
\end{bmatrix}\,,
\end{align}
the system (\ref{eq:dotq_ETSSC_v3}) and (\ref{eq:dhatG_ETSSC_v_3}) is simply reduced to
\begin{align}
\dfrac{dX(\bar{t})}{d\bar{t}}&=\dfrac{1}{\omega_{3}}\mathcal{F}\left(\bar{t},X,\dfrac{1}{\omega_{3}}\right)\,. \label{eq:dotX_event}
\end{align}
Note that (\ref{eq:dotX_event}) is characterized by a small parameter $1/\omega_{3}$ as well as a $T$-periodic function $\mathcal{F}\left(\bar{t},X,\dfrac{1}{\omega_{3}}\right)$ in $\bar{t}$ and, thereby, the averaging method can be performed on  $\mathcal{F}\left(\bar{t},X,\dfrac{1}{\omega_{3}}\right)$ at $\displaystyle \lim_{\omega\to \infty}\dfrac{1}{\omega_{3}}=0$, as shown in references \cite{K:2002,P:1979}. The averaging method allows for determining in what sense the behavior of a constructed average autonomous system approximates the behavior of the non-autonomous system (\ref{eq:dotX_event}). By employing the averaging technique to (\ref{eq:dotX_event}), we derive the following average system
\begin{align}
\dfrac{dX_{\rm{av}}(\bar{t})}{d\bar{t}}&=\dfrac{1}{\omega_{3}}\mathcal{F}_{\rm{av}}\left(X_{\rm{av}}\right) \,, \label{eq:dotXav_event_1} \\
\mathcal{F}_{\rm{av}}\left(X_{\rm{av}}\right)&=\dfrac{1}{T}\int_{0}^{T}\mathcal{F}\left(\gamma,X_{\rm{av}},0\right)d\gamma
\,,  \label{eq:mathcalFav_event}
\end{align}
where the terms with non-zero average values are
\begin{small}
\begin{align}
A_{\rm{av}}(\bar{t})&\mathbb{=} \frac{1}{T}\int_{0}^{T}A(\gamma)d\gamma=A \,, ~~  B_{\rm{av}}(\bar{t})\mathbb{=} \frac{1}{T}\int_{0}^{T}B(\gamma)d\gamma=B\,, \label{eq:Sav_event} \\
\!\!\!\mbox{\calligra H}~_{\rm{av}}(\bar{t})&\mathbb{=} \frac{1}{T}\int_{0}^{T}\!\!\!\!\mbox{\calligra H}~~(\gamma)d\gamma=I_{3} \,, ~~  \Delta_{\rm{av}}(\bar{t})\mathbb{=} \frac{1}{T}\int_{0}^{T}\Delta(\gamma)d\gamma=\bar{\Delta}\,, \label{eq:Mav_event} \\
\hat{\Delta}_{\rm{av}}(\bar{t})&\mathbb{=} \frac{1}{T}\int_{0}^{T}\hat{\Delta}(\gamma)d\gamma\mathbb{=}0\,, \label{eq:Deltaav_event}
\end{align}
\end{small}
and, consequently, 
\begin{align}
\mathcal{A}_{\rm{av}}(\bar{t})&=A \quad \mbox{ and } \quad
\mathcal{B}_{\rm{av}}(\bar{t})=B \,. \label{eq:dotHav_event}
\end{align}
Thus, using  (\ref{eq:Sav_event})--(\ref{eq:dotHav_event}), the average versions of (\ref{eq:dotq_ETSSC_v3}) and (\ref{eq:dhatG_ETSSC_v_3}) are, for all $\bar t\in [\bar t_k, \bar t_{k+1}):$
\begin{small}
\begin{align}
\frac{d\tilde{q}_{\rm{av}}}{d\bar{t}}(\bar{t})&\mathbb{=}\frac{1}{\omega_{3}}\left(A-BK\right)\tilde{q}_{\rm{av}}(\bar{t})-\frac{1}{\omega_{3}}BKe_{\rm{av}}(\bar{t})+\frac{1}{\omega_{3}}\bar{\Delta}\,, \label{eq:dotq_ETSSC_v4} \\
\frac{d\hat{G}_{\rm{av}}}{d\bar{t}}(\bar{t})&\mathbb{=}\frac{1}{\omega_{3}}A\tilde{q}_{\rm{av}}(\bar{t}) -\frac{1}{\omega_{3}}BK\hat{G}_{\rm{av}}(\bar{t})-\frac{1}{\omega_{3}}BKe_{\rm{av}}(\bar{t})\mathbb{+}\frac{1}{\omega_{3}}\bar{\Delta}\,. \label{eq:dhatG_ETSSC_v_4} 
\end{align}
\end{small}
Moreover, by using (\ref{eq:Sav_event})--(\ref{eq:dotHav_event}), the average version of (\ref{eq:hatG_ETSSC_v2}) and (\ref{eq:e_event}) are
\begin{align}
\hat{G}_{\rm{av}}(\bar{t})&= \tilde{q}_{\rm{av}}(\bar{t})\,, \label{eq:hatGav_event_1} \\
e_{\rm{av}}(\bar{t})&=\hat{G}_{\rm{av}}(\bar{t}_{k})-\hat{G}_{\rm{av}}(\bar{t})\,, \label{eq:Eav_event_1} 
\end{align}
and, consequently, the dynamics of $\tilde{q}_{\rm{av}}(\bar{t})$ and $\hat{G}_{\rm{av}}(\bar{t})$ are equivalent since, plugging (\ref{eq:hatGav_event_1}) into (\ref{eq:dhatG_ETSSC_v_4}), we obtain 
\begin{align}
\frac{d\hat{G}_{\rm{av}}}{d\bar{t}}(\bar{t})&=\frac{1}{\omega_{3}}(A-BK)\hat{G}_{\rm{av}}(\bar{t})-\frac{1}{\omega_{3}}BKe_{\rm{av}}(\bar{t})+\frac{1}{\omega_{3}}\bar{\Delta}\,, \label{eq:dhatG_ETSSC_v_5} 
\end{align}
preserving the  ISS property from $e_{\rm{av}}(\bar{t})$ and $\bar{\Delta}$ to $\hat{G}_{\rm{av}}(\bar{t})$ as well as $\tilde{q}_{\rm{av}}(\bar{t})$, see (\ref{eq:dotq_ETSSC_v4}).  Therefore, an average event-triggering can be introduced for the average system as well.

Defining the average version of $\Xi(\hat{G},e)$ as follows
\begin{small}
\begin{align}
\Xi(\hat{G}_{\rm{av}},e_{\rm{av}})&=\sigma\|\hat{G}_{\rm{av}}(\bar{t})\|\mathbb{-}\alpha(\|e_{\rm{av}}(\bar{t})\|\mathbb{+}a_{1}\omega_{3}\left|J_{2}(a_{3})\right|)	\,, \label{eq:Xi_event_2}
\end{align}
\end{small}
$\!\!$we construct the following average event-triggering mechanisms, according to the following definition.

\begin{definition}[\small{Average Static Triggering Condition}] 
\label{def:averageStaticEvent} Let $\Xi(\hat{G}_{\rm{av}},e_{\rm{av}})$ in (\ref{eq:Xi_event_2}) be the nonlinear mapping  and $K$  the control gain in (\ref{eq:U_event}). The event-triggered controller with average static-triggering condition in the new time scale $\bar{t}$ in (\ref{namadrugacomaltashoras}) consists of two components:
\begin{enumerate}
	\item A set of increasing sequence of time $I=\{\bar{t}_{0}\,, \bar{t}_{1}\,, \bar{t}_{2}\,,\ldots\}$ with $\bar{t}_{0}=0$ generated under the following rule:
		\begin{itemize}
			\item If $\left\{\bar{t} \in\mathbb{R}^{+}: \bar{t}>\bar{t}_{k} ~ \wedge ~ \Xi (\hat{G}_{\rm{av}},e_{\rm{av}}) < 0\right\} = \emptyset$, then the set of the times of the events is $I=\{\bar{t}_{0}\,, \bar{t}_{1}\,, \ldots, \bar{t}_{k}\}$.
			\item If $\left\{\bar{t} \in\mathbb{R}^{+}: \bar{t}>\bar{t}_{k} ~ \wedge ~ \Xi (\hat{G}_{\rm{av}},e_{\rm{av}}) < 0\right\} \neq \emptyset$, the next event time is given by
				\begin{align}
					\bar{t}_{k+1}&=\inf\left\{\bar{t} \in\mathbb{R}^{+}: \bar{t}>\bar{t}_{k} ~ \wedge ~ \Xi (\hat{G}_{\rm{av}},e_{\rm{av}}) <0 \right\}\,, \label{eq:tk+1_event_av}
				\end{align}
				\!\!\!\!\!\!\!\!\!\!\!\!which is the average static event-triggering mechanism.
		\end{itemize}
		\item An average feedback control action updated at the triggering instants:
		\begin{align}
			u_{k}^{\rm av}=-K\hat{G}_{\rm av}(\bar{t}_{k}) \,, \label{eq:U_MD4}
		\end{align}
		for all $\bar{t} \in \lbrack \bar{t}_{k}\,, \bar{t}_{k+1}\phantom{(}\!\!)$, $k\in \mathbb{N}$.
\end{enumerate}
\end{definition}

The event-triggering mechanism discussed above guarantees the asymptotic stabilization of $\hat{G}_{\rm{av}}(\bar{t})$ and $\tilde{q}_{\rm{av}}(\bar{t})$. Consequently, both $\hat{G}_{\rm{av}}(\bar{t})$ and $\tilde{q}_{\rm{av}}(\bar{t})$ converge to a neighborhood of origin according to the averaging theory \cite{K:2002}.

\section{STABILITY ANALYSIS}


Theorem~\ref{thm:NETESC_4} states the local asymptotic stability of the average closed-loop system and the exponential convergence to a neighborhood of the extremum point by means of our static ET-SSC strategy.

\begin{theorem} \label{thm:NETESC_4}
Consider the closed-loop average dynamics of the gradient estimate (\ref{eq:dhatG_ETSSC_v_5}), the average error vector (\ref{eq:Eav_event_1})  and the average \textbf{static} event-triggering mechanism given by \textbf{Definition \ref{def:averageStaticEvent}}. Under  Assumptions \ref{assumption1}--\ref{assumption2} and considering the quadratic mapping $\Xi(\hat{G}_{\rm{av}},e_{\rm{av}})$ given by (\ref{eq:Xi_event_2}), for  $\omega_{3}>0$, defined in (\ref{eq:w1_w2}), sufficiently large, the equilibrium $\hat{G}_{\rm{av}}(t)=0$ is locally exponentially stable and $\tilde{q}_{\rm{av}}(t)$ converges exponentially to zero such that 
\begin{align}
&\|q(t) \mathbb{-} q^\ast\|  \leq  \nonumber \\
&\exp\left(\!\!\mathbb{-}\frac{1}{2}\frac{\lambda_{\min}(Q)}{\lambda_{\max}(P)}(1\mathbb{-}\sigma)t\right)\sqrt{\frac{\lambda_{\max}(P)}{\lambda_{\min}(P)}}\|q(0) \mathbb{-} q^\ast\| \nonumber \\
&+\mathcal{O}\left(a+\frac{1}{\omega}\right) \,, 
\label{eq:normQ_thm}
\end{align}
where $a=\frac{3}{2}a_{3}$, and $\omega = \omega_{3}$. 
 In addition, there exists a lower bound  $\tau^{\ast}$ for the inter-event interval $t_{k+1}-t_{k}$  for all $k \in \mathbb{N}$ precluding the Zeno behavior.
\end{theorem}
\textit{Proof:} See Appendix. \hfill $\square$ \\ 

\section{SIMULATION RESULTS}

In the following simulations, we set the parameters of the
stationary source as $[x^*\,, y^*]=[10\,, 5]$ meters, whereas $\theta^* =\pi/6$ rad, and $Q^*= 7$. The parameters of the proposed ET-SSC are
chosen as $\omega_1=\omega_2=10$ rad/s, $\omega_3=20$ rad/s,  $a_1=a_2=a_3=0.5$, $K=\begin{bmatrix} 4.3822 & 4.3822 & 0.1437 \\ -9.4326 & 9.4326 & 4\end{bmatrix}$, $\sigma=0.5$, and $\alpha=0.195$. The
starting position of the autonomous vehicle is  $[x(0)\,, y(0)]=[12.5\,, 7.5]$ meters, with initial orientation $\theta(0)=\pi/3$ rad.

\begin{figure}[h!]
	\centering
    	\subfigure[Nonlinear map output, $Q(t)$. \label{fig:Q}]{\includegraphics[width=4.26cm]{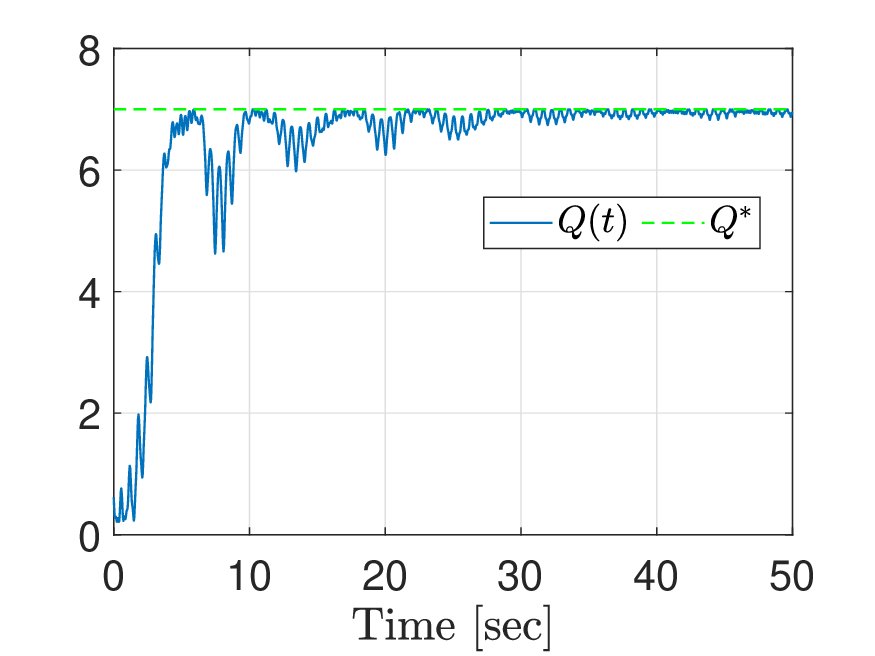}}
	\subfigure[Nonlinear map input, $q(t)$. \label{fig:q}]{\includegraphics[width=4.26cm]{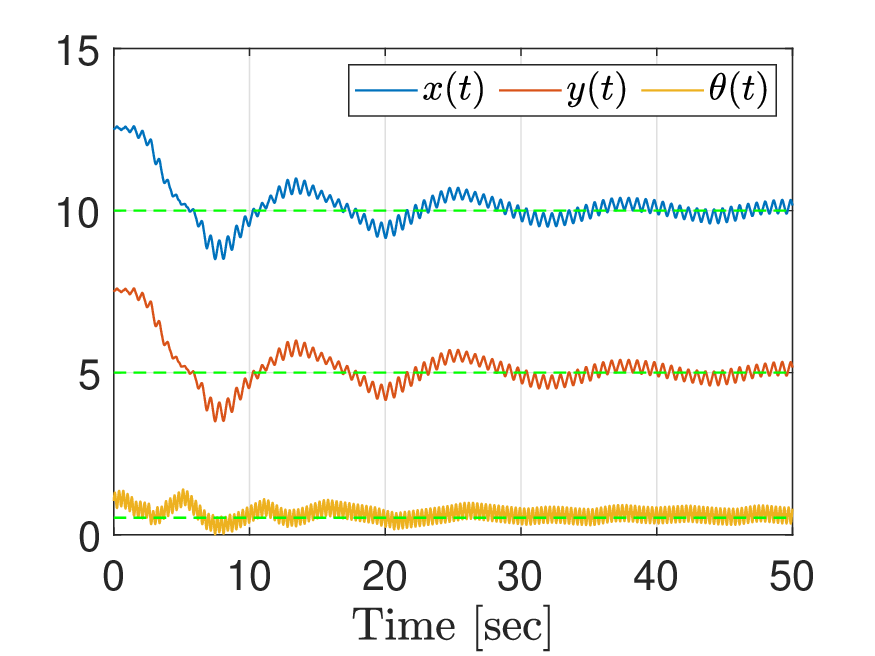}}
		\\
        		\subfigure[Gradient estimate, $\hat{G}(t_{k})$. \label{fig:G}]{\includegraphics[width=4.26cm]{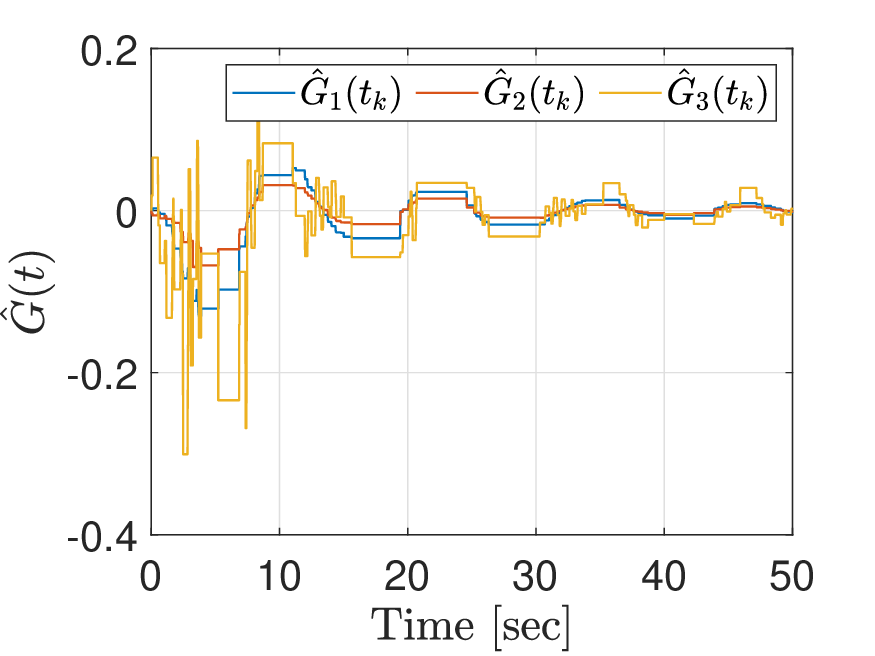}}
		\subfigure[Control signal, $u(t)$. \label{fig:u}]{\includegraphics[width=4.26cm]{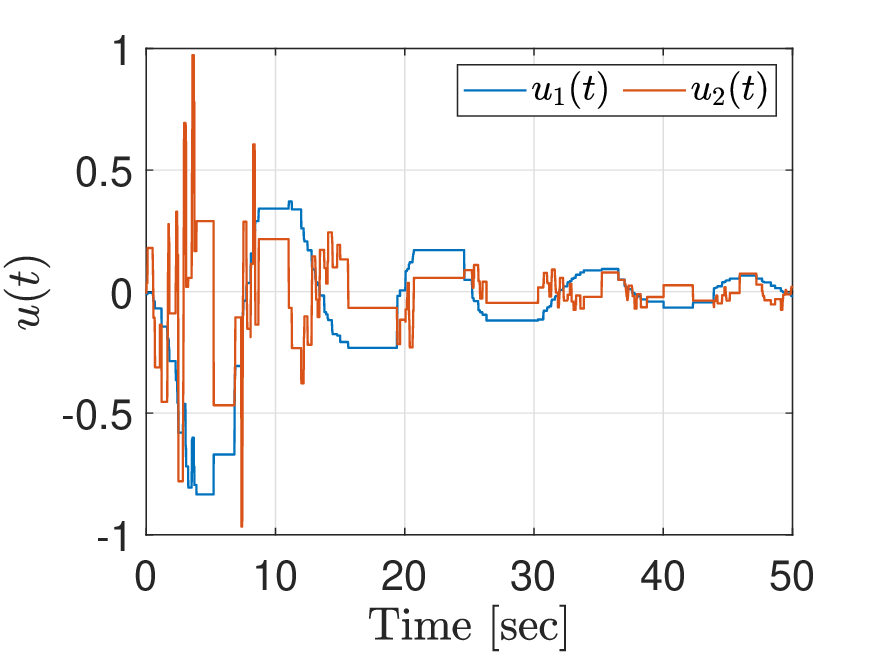}}
\end{figure}

The output of the unknown signal $Q(x\,,y\,,\theta)$ is shown in
Fig.~\ref{fig:Q}, while the trajectories and orientation of the vehicle converging to the maximizer $[x^{\ast}, y^{\ast}, \theta^{\ast}]$ are shown Fig.~\ref{fig:q}. In Fig.~\ref{fig:G}, the gradient estimates are are being updated at triggering times ensuring their convergence to zero. Figure~\ref{fig:u} highlights the aperiodic behavior of how often the control signal $u(t)$ is updated, showing that the proposed ET approach achieves the optimal value with fewer control signal updates than the continuous SSC or even a sampled-data version would be able to perform. 

\section{CONCLUSIONS}
This paper presented an event-triggered source seeking control (ET-SSC) strategy for autonomous vehicles modeled as nonholonomic unicycles. By introducing static-triggering conditions, we reduced the frequency of control updates, leading to a resource-efficient design without compromising performance. 
Through a combination of time-scaling and Lyapunov methods with averaging theory, we proved that the proposed approach guaranteed exponential stability of the average system, resulting in practical asymptotic convergence to a small neighborhood of the source. Additionally, we established a minimum dwell time to eliminate Zeno behavior, guaranteeing reliable implementation. Compared to conventional continuous-time source seeking algorithms, ET-SSC achieved a similar level of optimization while significantly reducing actuation demands. A key advantage of the proposed method was its ability to accommodate arbitrarily large inter-sampling intervals, addressing key limitations of existing sampled-data approaches. 
Numerical simulations corroborated our findings, illustrating the effectiveness and practicality of the control strategy.

Future investigation lies in the design and analysis of different control problems with event-triggered implementation, as considered in the following references \cite{paper1,paper2,paper3,paper4,paper5,paper6,paper7,paper8,paper9,paper10,paper11,paper12,paper13,paper14,paper15,paper16,paper17,paper18,paper19,paper20}.

\begin{small}
\section*{Appendix:\\ Proof of Theorem~\ref{thm:NETESC_4}} \label{appendix:dynamic}

The proof of the theorem is again divided into two parts: stability analysis and avoidance of Zeno behavior.

\subsection{Stability Analysis}

Consider the following Lyapunov candidate for the   average system:
\begin{align}
V_{\rm{av}}(\bar{t})=\hat{G}^{T}_{\rm{av}}(\bar{t})P\hat{G}_{\rm{av}}(\bar{t})\,,\quad P^{T}=P>0\,. \label{eq:lyapunov_dynamicETC_1_pf4}
\end{align}
The Rayleigh-Ritz inequality writes:
\begin{align}
\lambda_{\min}(P)\|\hat{G}_{\rm{av}}(\bar{t})\|^2 \leq \hat{G}^{T}_{\rm{av}}(\bar{t})P\hat{G}_{\rm{av}}(\bar{t}) \leq \lambda_{\max}(P)\|\hat{G}_{\rm{av}}(\bar{t})\|^2 \,. \label{ineq:RRI_dynamic_pf4}
\end{align}
The time-derivative of (\ref{eq:lyapunov_dynamicETC_1_pf4}) is given by
\begin{align}
&\frac{d{V}_{\rm{av}}(\bar{t})}{d\bar{t}}=\frac{d\hat{G}_{\rm{av}}^{T}(\bar{t})}{d\bar{t}}P\hat{G}_{\rm{av}}(\bar{t})\mathbb{+}\hat{G}^{T}_{\rm{av}}(\bar{t})P\frac{d\hat{G}_{\rm{av}}(\bar{t})}{d\bar{t}}\,,  \label{eq:dotLyapunov_dynamicETC_1_pf4}
\end{align} 
which, by using equations (\ref{eq:dhatG_ETSSC_v_5}) and (\ref{eq:LyapEq}), can be rewritten as
\begin{align}
\frac{d{V}_{\rm{av}}(\bar{t})}{d\bar{t}}&\mathbb{=}\mathbb{-}\frac{1}{\omega_{3}}\hat{G}_{\rm{av}}^{T}(\bar{t})Q\hat{G}_{\rm{av}}(\bar{t})\mathbb{+}\frac{2}{\omega_{3}}\hat{G}_{\rm{av}}^{T}(\bar{t})P(A\mathbb{-}BK)e_{\rm{av}}(\bar{t}) \nonumber \\
&\quad\mathbb{+}\frac{2}{\omega_{3}}\hat{G}_{\rm{av}}^{T}(\bar{t})P(A\mathbb{-}BK)\bar{\Delta}\,. \label{eq:dotLyapunov_dynamicETC_2_pf4}
\end{align}
By using (\ref{eq:barDelta}), an upper bound for the norm of the average disturbance $\bar{\Delta}$ can be obtained using only design parameters from the extremum seeking scheme, along with the frequency and amplitudes of the perturbation signals. This upper bound satisfy 
\begin{align}
\|\bar{\Delta}\|
    & \leq a_{1}\omega_{3}\left|J_{2}(a_{3})\right| \,. \label{eq:barDelta_v2}
\end{align}
Therefore, by using the upper bound (\ref{eq:barDelta_v2}), the lower bound $\lambda_{\min}(Q)\|\hat{G}_{\rm{av}}(\bar{t})\|^{2}\leq \hat{G}_{\rm{av}}^{T}(\bar{t})Q\hat{G}_{\rm{av}}(\bar{t})$ and inequality (\ref{eq:alpha}), equation (\ref{eq:dotLyapunov_dynamicETC_2_pf4}) is upper bounded by
\begin{small}
\begin{align}
&\frac{d{V}_{\rm{av}}(\bar{t})}{d\bar{t}}\mathbb{\leq}\mathbb{-}\frac{\lambda_{\min}(Q)}{\omega_{3}}\|\hat{G}_{\rm{av}}(\bar{t})\|^{2}\mathbb{+} \nonumber \\
&\mathbb{+}\frac{2}{\omega_{3}}\|\hat{G}_{\rm{av}}(\bar{t})\|\|P(A\mathbb{-}BK)\|(\|e_{\rm{av}}(\bar{t})\|\mathbb{+}a_{1}\omega_{3}\left|J_{2}(a_{3})\right|) \nonumber \\
&=\mathbb{-}\frac{\lambda_{\min}(Q)}{\omega_{3}}\|\hat{G}_{\rm{av}}(\bar{t})\|\times  \nonumber \\
&\left(\|\hat{G}_{\rm{av}}(\bar{t})\|\mathbb{-}\frac{2\|P(A\mathbb{-}BK)\|}{\lambda_{\min}(Q)}(\|e_{\rm{av}}(\bar{t})\|\mathbb{+}a_{1}\omega_{3}\left|J_{2}(a_{3})\right|)\right)\\
&\mathbb{\leq}\mathbb{-}\frac{\lambda_{\min}(Q)}{\omega_{3}}\|\hat{G}_{\rm{av}}(\bar{t})\|\times  \nonumber \\
&\quad\left(\|\hat{G}_{\rm{av}}(\bar{t})\|\mathbb{-}\alpha(\|e_{\rm{av}}(\bar{t})\|\mathbb{+}a_{1}\omega_{3}\left|J_{2}(a_{3})\right|)\right)\,. \label{eq:dotV_ETSSC_v4}
\end{align}
\end{small}

In the proposed event-triggered mechanism, the update law is (\ref{eq:tk+1_event_av}) and $\Xi(\hat{G}_{\rm{av}},e_{\rm{av}})$ is given by (\ref{eq:Xi_event_2}). The signal $u_{\text{av}}(t)$ is held constant between two consecutive events, {\it i.e.}, while $\Xi(\hat{G}_{\rm{av}},e_{\rm{av}})\geq 0$, consequently
\begin{align}
\sigma\|\hat{G}_{\rm{av}}(\bar{t})\|\mathbb{-}\alpha(\|e_{\rm{av}}(\bar{t})\|\mathbb{+}a_{1}\omega_{3}\left|J_{2}(a_{3})\right|)\geq 0		\,.
\end{align}
Therefore, the following upper bound is verified 
\begin{align}
\alpha(\|e_{\rm{av}}(\bar{t})\|\mathbb{+}a_{1}\omega_{3}\left|J_{2}(a_{3})\right|)&\leq \sigma\|\hat{G}_{\rm{av}}(\bar{t})\|	\,. \label{eq:eAv_upperBound}
\end{align}
Then, by using (\ref{eq:eAv_upperBound}), inequality (\ref{eq:dotV_ETSSC_v4}) is upper bounded as
\begin{align}
\frac{d{V}_{\rm{av}}(\bar{t})}{d\bar{t}}&\mathbb{\leq}\mathbb{-}\frac{\lambda_{\min}(Q)}{\omega_{3}}(1-\sigma)\|\hat{G}_{\rm{av}}(\bar{t})\|^{2}\,, \quad \forall \bar{t}\in \lbrack \bar{t}_{k},\bar{t}_{k+1}\phantom{(}\!\!) \,. \label{eq:dotV_ETSSC_v5}
\end{align}

Now,  using  (\ref{ineq:RRI_dynamic_pf4}), inequality (\ref{eq:dotV_ETSSC_v5}) can be upper bounded as follows
\begin{align}
\frac{d{V}_{\rm{av}}(\bar{t})}{d\bar{t}}&\mathbb{\leq}\mathbb{-}\frac{1}{\omega_{3}}\frac{\lambda_{\min}(Q)}{\lambda_{\max}(P)}(1-\sigma)V_{\rm{av}}(\bar{t})\,, \quad \forall \bar{t}\in \lbrack \bar{t}_{k},\bar{t}_{k+1}\phantom{(}\!\!) \,. \label{eq:dotV_ETSSC_v6}
\end{align}

Then, invoking the Comparison Principle \cite[Comparison Lemma]{K:2002}, an upper bound for $\bar{V}_{\rm{av}}(\bar{t})$ for $V_{\rm{av}}(\bar{t})$ 
\begin{align}
V_{\rm{av}}(\bar{t}) \leq \bar{V}_{\rm{av}}(\bar{t})\,, \quad \forall \bar{t}\in \lbrack \bar{t}_{k},\bar{t}_{k+1}\phantom{(}\!\!)\,, \label{ineq:VavBarVav_dynamic_pf4}
\end{align} 
is given by the solution of the dynamics
\begin{align}
\frac{d\bar{V}_{\rm{av}}(\bar{t})}{d\bar{t}}&\mathbb{=}\mathbb{-}\frac{1}{\omega_{3}}\frac{\lambda_{\min}(Q)}{\lambda_{\max}(P)}(1-\sigma)\bar{V}_{\rm{av}}(\bar{t}),~\bar{V}_{\rm{av}}(\bar{t}_{k})\mathbb{=}V_{\rm{av}}(\bar{t}_{k}). \label{eq:barVav_0_dynamic_pf4}
\end{align}
Hence, $\forall \bar{t}\in \lbrack \bar{t}_{k},\bar{t}_{k+1}\phantom{(}\!\!)$, one has:
\begin{align}
\bar{V}_{\rm{av}}(\bar{t})=\exp\left(\mathbb{-}\frac{1}{\omega_{3}}\frac{\lambda_{\min}(Q)}{\lambda_{\max}(P)}(1-\sigma)\bar{t}\right)\bar{V}_{\rm{av}}(\bar{t}_{k})\,. \label{eq:barVav_dynamic_pf4}
\end{align}
Using (\ref{eq:barVav_0_dynamic_pf4}) and (\ref{eq:barVav_dynamic_pf4}), the inequality (\ref{ineq:VavBarVav_dynamic_pf4}) is rewritten as 
\begin{align}
V_{\rm{av}}(\bar{t})\leq\exp\left(\mathbb{-}\frac{1}{\omega_{3}}\frac{\lambda_{\min}(Q)}{\lambda_{\max}(P)}(1-\sigma)\bar{t}\right)V_{\rm{av}}(\bar{t}_{k})\,. \label{ineq:Vav_dynamic_pf4}
\end{align}
By defining, $\bar{t}_{k}^{+}$ and $\bar{t}_{k}^{-}$ as the right and left limits of $\bar{t}=\bar{t}_{k}$, respectively, it easy to verify that {\small$$V_{\rm{av}}(\bar{t}_{k+1}^{-})\leq \exp\left(\mathbb{-}\frac{1}{\omega_{3}}\frac{\lambda_{\min}(Q)}{\lambda_{\max}(P)}(1-\sigma)(\bar{t}_{k+1}^{-}-\bar{t}_{k}^{+})\right)V_{\rm{av}}(\bar{t}_{k}^{+}).$$} Since $V_{\rm{av}}(\bar{t})$ is continuous, $V_{\rm{av}}(\bar{t}_{k+1}^{-})=V_{\rm{av}}(\bar{t}_{k+1})$ and $V_{\rm{av}}(\bar{t}_{k}^{+})=V_{\rm{av}}(\bar{t}_{k})$, 
and therefore,
\begin{align}
    V_{\rm{av}}(\bar{t}_{k+1})\mathbb{\leq} \exp\left(\mathbb{-}\frac{1}{\omega_{3}}\frac{\lambda_{\min}(Q)}{\lambda_{\max}(P)}(1\mathbb{-}\sigma)(\bar{t}_{k+1}\mathbb{-}\bar{t}_{k})\right)V_{\rm{av}}(\bar{t}_{k}). \label{METES_eq:mmd_2_s}
\end{align}
To handle the discontinuities in the Lyapunov function at each triggering event, we have defined the Lyapunov function piecewise over intervals between triggering times. Then, we analyze the Lyapunov function's value just before and after each triggering event to manage the jumps. As done in \cite{s5}, by establishing a recursive analysis that shows how the Lyapunov function decreases across these intervals, we demonstrate that the Lyapunov function's value after each jump is consistently lower than before. This ensures an overall exponential decrease in the Lyapunov function over all time, as follows. 

%
Hence, for any $\bar{t}\geq 0$ in $ \bar{t}\in \lbrack \bar{t}_{k},\bar{t}_{k+1}\phantom{(}\!\!)$, $k \in \mathbb{N}$, one has 
\begin{align}
    &V_{\rm{av}}(\bar{t})\leq \exp\left(\mathbb{-}\frac{1}{\omega_{3}}\frac{\lambda_{\min}(Q)}{\lambda_{\max}(P)}(1-\sigma)(\bar{t}-\bar{t}_{k})\right) V_{\rm{av}}(\bar{t}_{k}) \nonumber \\
    &\leq \exp\left(\mathbb{-}\frac{1}{\omega_{3}}\frac{\lambda_{\min}(Q)}{\lambda_{\max}(P)}(1-\sigma)(\bar{t}-\bar{t}_{k})\right)  \nonumber \\
		& \times \exp\left(\mathbb{-}\frac{1}{\omega_{3}}\frac{\lambda_{\min}(Q)}{\lambda_{\max}(P)}(1-\sigma)(\bar{t}_{k}-\bar{t}_{k-1})\right)V_{\rm{av}}(\bar{t}_{k-1}) \nonumber \\
    &\leq \ldots \leq \nonumber \\
    &\leq \exp\left(\mathbb{-}\frac{1}{\omega_{3}}\frac{\lambda_{\min}(Q)}{\lambda_{\max}(P)}(1-\sigma)(\bar{t}\mathbb{-}\bar{t}_{k})\right) \nonumber \\
		& \times \prod_{i=1}^{i=k}\exp\left(\mathbb{-}\frac{1}{\omega_{3}}\frac{\lambda_{\min}(Q)}{\lambda_{\max}(P)}(1-\sigma)(\bar{t}_{i}-\bar{t}_{i-1})\right)V_{\rm{av}}(\bar{t}_{i-1}) \nonumber \\
    &=\exp\left(\mathbb{-}\frac{1}{\omega_{3}}\frac{\lambda_{\min}(Q)}{\lambda_{\max}(P)}(1-\sigma)\bar{t}\right) V_{\rm{av}}(0)\,. \label{METES_eq:barVav_dynamic_pf4}
\end{align}
From (\ref{eq:lyapunov_dynamicETC_1_pf4}), it follows 
\begin{align}
\hat{G}^{T}_{\rm{av}}(\bar{t})P\hat{G}_{\rm{av}}(\bar{t})\leq V_{\rm{av}}(\bar{t})\,. \label{ineq:hatGav_Vav_dynamic_1_pf4}
\end{align}
Consequently, combining  (\ref{METES_eq:barVav_dynamic_pf4}) and (\ref{ineq:hatGav_Vav_dynamic_1_pf4}), one gets
\begin{align}
&\hat{G}^{T}_{\rm{av}}(\bar{t})P\hat{G}_{\rm{av}}(\bar{t})\leq \exp\left(\mathbb{-}\frac{1}{\omega_{3}}\frac{\lambda_{\min}(Q)}{\lambda_{\max}(P)}(1-\sigma)\bar{t}\right) V_{\rm{av}}(0)\,, \nonumber \\
&= \exp\left(\mathbb{-}\frac{1}{\omega_{3}}\frac{\lambda_{\min}(Q)}{\lambda_{\max}(P)}(1-\sigma)\bar{t}\right)\hat{G}^{T}_{\rm{av}}(0)P\hat{G}_{\rm{av}}(0)\,. \label{ineq:hatGav_Vav_dynamic_2_pf4}
\end{align}
Therefore, from (\ref{ineq:RRI_dynamic_pf4}), one gets
\begin{small}
\begin{align}
\|\hat{G}_{\rm{av}}(\bar{t})\|^2 & \!\leq \! \exp\left(\mathbb{-}\frac{1}{\omega_{3}}\frac{\lambda_{\min}(Q)}{\lambda_{\max}(P)}(1-\sigma)\bar{t}\right)\frac{\lambda_{\max}(P)}{\lambda_{\min}(P)}\|\hat{G}_{\rm{av}}(0)\|^2\!\!. \label{ineq:hatGav_Vav_dynamic_3_pf4}
\end{align}
\end{small}
and
\begin{small}
\begin{align}
\|\hat{G}_{\rm{av}}(\bar{t})\| & \!\leq \! \exp\left(\!\!\mathbb{-}\frac{1}{2\omega_{3}}\frac{\lambda_{\min}(Q)}{\lambda_{\max}(P)}(1\mathbb{-}\sigma)\bar{t}\right)\sqrt{\frac{\lambda_{\max}(P)}{\lambda_{\min}(P)}}\|\hat{G}_{\rm{av}}(0)\|\,. \label{ineq:hatGav_Vav_dynamic_7_pf4}
\end{align}
\end{small}
Although the analysis has been focused on the convergence of $\hat{G}_{\rm{av}}(\bar{t})$ and, consequently, $\hat{G}(t)$, the obtained results through (\ref{ineq:hatGav_Vav_dynamic_7_pf4}) can be easily extended to the variables $\tilde{q}_{\rm{av}}(\bar{t})$ and $\tilde{q}(t)$ since, from (\ref{eq:hatGav_event_1}), $\hat{G}_{\rm{av}}(\bar{t})=\tilde{q}_{\rm{av}}(\bar{t})$, therefore, inequality (\ref{ineq:hatGav_Vav_dynamic_7_pf4}) can be rewritten as
\begin{small}
\begin{align}
\|\tilde{q}_{\rm{av}}(\bar{t})\| & \!\leq \! \exp\left(\!\!\mathbb{-}\frac{1}{2\omega_{3}}\frac{\lambda_{\min}(Q)}{\lambda_{\max}(P)}(1\mathbb{-}\sigma)\bar{t}\right)\sqrt{\frac{\lambda_{\max}(P)}{\lambda_{\min}(P)}}\|\tilde{q}_{\rm{av}}(0)\|\,. \label{ineq:tildeQav_v1}
\end{align}
\end{small}

Since the right-hand side of the differential equation~(\ref{eq:dotq_ETSSC_v3}) is discontinuous due to the control action, $T$-periodic in $t$ due to the probing signals, and satisfies the Lipschitz condition, and given that $\tilde{q}_{\rm{av}}(\bar{t})$ is asymptotically stable as established by~(\ref{ineq:tildeQav_v1}), we can invoke~\cite[Theorem~2]{P:1979}, such that  
\begin{align}
\|\tilde{q}(\bar{t})-\tilde{q}_{\rm{av}}(\bar{t})\|\leq\mathcal{O}\left(\frac{1}{\omega}\right)\,. \label{eq:plotnikov}
\end{align}

Now, adding and subtracting $\tilde{q}_{\rm{av}}(\bar{t})$ in the right-hand side of (\ref{eq:q_ETSSC_v2}), in time scale $\bar{t}$m one has
\begin{align}
q(\bar{t}) - q^\ast = \tilde{q}_{\rm{av}}(\bar{t})+\tilde{q}(\bar{t})-\tilde{q}_{\rm{av}}(\bar{t}) + S(\bar{t})\,, \label{eq:q_ETSSC_v3cu}
\end{align}  
whose norm can be upper bounded, by using the triangle inequality \cite{A:1957}, such that
\begin{align}
\|q(\bar{t}) - q^\ast\| \leq \|\tilde{q}_{\rm{av}}(\bar{t})\|+\|\tilde{q}(\bar{t})-\tilde{q}_{\rm{av}}(\bar{t})\| + \|S(\bar{t})\|\,.\label{eq:q_ETSSC_v3cu2}
\end{align} 
Since sine functions in (\ref{eq:S_v1}) are uniformly bounded, it is possible to derive a uniform upper bound for the Euclidean norm of $S(\bar{t})$ such that
\begin{align}
    \|S(\bar{t})\|\leq \frac{1}{2}\sqrt{a_{1}^{2}+a_{2}^{2}+a_{3}^{2}}=a \,.\label{eq:S_v2}
\end{align}
Thus, by using (\ref{ineq:tildeQav_v1}), (\ref{eq:plotnikov}) and (\ref{eq:S_v2}), inequality (\ref{eq:q_ETSSC_v3cu2}) is upper bounded by 
\begin{align}
&\|q(\bar{t}) \mathbb{-} q^\ast\|  \leq  \nonumber \\
& \exp\left(\!\!\mathbb{-}\frac{1}{2\omega_{3}}\frac{\lambda_{\min}(Q)}{\lambda_{\max}(P)}(1\mathbb{-}\sigma)\bar{t}\right)\sqrt{\frac{\lambda_{\max}(P)}{\lambda_{\min}(P)}}\|q(0) \mathbb{-} q^\ast\| \nonumber \\
&+\mathcal{O}\left(a+\frac{1}{\omega}\right)\,.\label{eq:q_ETSSC_v4}
\end{align} 
Therefore, by using (\ref{eq:w1_w2}) and the original system time scale $t=\frac{\bar{t}}{\omega_{3}}$, inequality (\ref{eq:normQ_thm}) is verified.

\subsection{Avoidance of Zeno Behavior}

Since the average closed-loop system consists of (\ref{eq:dhatG_ETSSC_v_5}), with the event-triggering mechanism  (\ref{eq:tk+1_event_av}) and the average control law (\ref{eq:U_MD4}), we can conclude that $\alpha(\|e_{\rm{av}}(\bar{t})\|\mathbb{+}8(a_{1}\mathbb{+}a_{2})\omega_{3}\left|J_{2}(a_{3})\right|)| \leq \sigma\|\hat{G}_{\text{av}}(\bar{t})\|$, resulting in  
\begin{small}
\begin{align}
\sigma\|\hat{G}_{\text{av}}(\bar{t})\|^{2}-\alpha(\|e_{\rm{av}}(\bar{t})\|\mathbb{+}a_{1}\omega_{3}\left|J_{2}(a_{3})\right|)\|\hat{G}_{\text{av}}(\bar{t})\|\geq 0\,. \label{ineq:interEvents_1_static}
\end{align}
\end{small}
By defining $\|\bar{e}_{\text{av}}(\bar{t})\|:=\alpha(\|e_{\rm{av}}(\bar{t})\|\mathbb{+}a_{1}\omega_{3}\left|J_{2}(a_{3})\right|)$ and using the Peter-Paul inequality \cite{W:1971}, $cd\leq \frac{c^2}{2\epsilon}+\frac{\epsilon d^2}{2}$ for all $c,d,\epsilon>0$, with $c=\|\bar{e}_{\rm{av}}(\bar{t})\|$, $d=\|\hat{G}_{\rm{av}}(\bar{t})\|$ and $\epsilon=\sigma$, the inequality (\ref{ineq:interEvents_1_static}) is lower bounded by
\begin{align}
&\sigma \|\hat{G}_{\text{av}}(\bar{t})\|^{2}-\|\bar{e}_{\text{av}}(\bar{t})\|\|\hat{G}_{\text{av}}(\bar{t})\|\geq \nonumber \\
&n\|\hat{G}_{\rm{av}}(\bar{t})\|^{2}-m\|\bar{e}_{\rm{av}}(\bar{t})\|^2\,,\label{ineq:interEvents_2_static_pf2}
\end{align}
where  $n=\frac{\sigma}{2}$ and $m=\frac{1}{2\sigma}$. 
In \cite{G:2014}, it is shown that a lower bound for the inter-execution interval is given by the time duration it takes for the function
\begin{align}
\phi_{\rm{av}}(\bar{t})=\sqrt{\frac{m}{n}}\frac{\|\bar{e}_{\rm{av}}(\bar{t})\|}{\|\hat{G}_{\rm{av}}(\bar{t})\|} \label{eq:phi_1_static_pf2}
\end{align}
to go from 0 to 1. The time-derivative of (\ref{eq:phi_1_static_pf2}) satisfies 
\begin{align}
\omega\frac{d\phi_{\text{av}}(\bar{t})}{d\bar{t}}&\leq\frac{\|A-BK\|+\|BK\|}{\omega}\sqrt{\frac{n}{m}}\left(\sqrt{\frac{m}{n}}+\phi_{\text{av}}(\bar{t})\right)^{2}\,. \label{eq:dotPhi_3_dynamic_pf4}
\end{align}
Then, invoking the Comparison Principle \cite[Comparison Lemma]{K:2002} an upper bound $\tilde{\phi}_{\text{av}}(\bar{t})$ for $\phi_{\text{av}}(\bar{t})$
\begin{align}
\phi_{\text{av}}(\bar{t})\leq \tilde{\phi}_{\text{av}}(\bar{t}) \,, \quad \phi_{\text{av}}(0)= \tilde{\phi}_{\text{av}}(0)=0 \,, \quad  \forall \bar{t}\in \lbrack \bar{t}_{k},\bar{t}_{k+1}\phantom{(}\!\!) \,, \label{eq:tildePhi_v2}
\end{align}
 is given by the solution of the equation
\begin{align}
\frac{d\tilde{\phi}_{\text{av}}(\bar{t})}{d\bar{t}}&=\frac{\|A-BK\|+\|BK\|}{\omega}\sqrt{\frac{n}{m}}\left(\sqrt{\frac{m}{n}}+\tilde{\phi}_{\text{av}}(\bar{t})\right)^{2}\,. \label{eq:dotTildePhi_v2}
\end{align}
The solution of (\ref{eq:dotTildePhi_v2}), with the initial condition $\tilde{\phi}_{\text{av}}(0) = 0$, is 
\begin{align}
\tilde{\phi}_{\text{av}}(\bar{t}) = \frac{\sqrt{\frac{m}{n}}}{1 - \frac{\|A-BK\|+\|BK\|}{\omega} \frac{n}{m} \bar{t}} - \sqrt{\frac{m}{n}}. \label{eq:phi-bart_v2}
\end{align}
This solution describes the evolution of $\tilde{\phi}_{\text{av}}(\bar{t})$ over time, starting from zero at $\bar{t} = 0$. Since $\phi(t)=\sqrt{\frac{m}{n}}\frac{\|\bar{e}(t)\|}{\|\hat{G}(t)\|}$, with $\bar{e}(t)$ and $\hat{G}(t)$ being $T$-periodic in $t$, and its average version, given by (\ref{eq:phi_1_static_pf2}), is upper bounded by (\ref{eq:phi-bart_v2}), by invoking \cite[Theorem 2]{P:1979},
\begin{align}
|\phi(t)-\tilde{\phi}_{\text{av}}(t)|\leq\mathcal{O}\left(\frac{1}{\omega}\right)\,.
\end{align}
By using the Triangle inequality \cite{A:1957}, one has
\begin{align}
\phi(t)&\leq\ \phi_{\text{av}}(t)+\mathcal{O}\left(\frac{1}{\omega}\right)\leq \tilde{\phi}_{\text{av}}(t)+\mathcal{O}\left(\frac{1}{\omega}\right)\! \nonumber \\
&=\frac{\sqrt{\frac{m}{n}}}{1 - (\|A-BK\|+\|BK\|) \frac{n}{m} t} - \sqrt{\frac{m}{n}} +\mathcal{O}\left(\frac{1}{\omega}\right). \label{eq:phi_t}
\end{align}
Now, defining 
\begin{align}
\hat{\phi}(t):=\frac{\sqrt{\frac{m}{n}}}{1 - (\|A-BK\|+\|BK\|) \frac{n}{m} t} - \sqrt{\frac{m}{n}} +\mathcal{O}\left(\frac{1}{\omega}\right)\,, \label{eq:hatPhi_t}
\end{align}
a lower bound for the inter-execution time of original system is given by the time it takes for the function (\ref{eq:hatPhi_t}) go to 0 to 1, this is at least
\begin{align}
\tau^{\ast}&=\frac{1}{\|A-BK\|+\|BK\|}\frac{m}{n}\frac{1-\mathcal{O}(1/\omega)}{1+\sqrt{m/n}-\mathcal{O}(1/\omega)}\,,
\end{align} 
and the Zeno behavior is avoided in the original system.
\end{small}


\begin{thebibliography}{99}

\bibitem{s9} 
K.~E.~{\r A}arz{\' e}n, 
\newblock A simple event-based PID controller.
\newblock {\em IFAC World Congress}, 32:423--428, 1999.

\bibitem{APDNH:2018} M. Abdelrahim, R. Postoyan, J. Daafouz, D. Ne{\u s}i{\' c}, and M. Heemels, "Co-design of output feedback laws and event-triggering conditions for the L2-stabilization of linear systems," {\it Automatica}, vol. 87, pp. 337--344, 2018.

\bibitem{APDN:2016}
M.~Abdelrahim, R.~Postoyan, J.~Daafouz, and D.~Ne{\u s}i{\' c}, 
\newblock Stabilization of nonlinear systems using event-triggered output feedback controllers.
\newblock {\em IEEE Trans. Automat. Contr.}, 61:2682--2687, 2016.


\bibitem{AS:1964} M. Abramowitz and I. A. Stegun, ``Handbook of mathematical functions with formulas, graphs, and mathematical tables,'' in \textit{Dover Books on Mathematics}, Revised ed., New York: Dover Publications, 1964.

\bibitem{AB:2011} 
T. Alpcan and T. Başar, Network Security: A Decision and Game Theoretic Approach, Cambridge University Press, 2011.


\bibitem{paper6}
N.~O. Aminde, T.~R. Oliveira, and L.~Hsu, Global output-feedback extremum
  seeking control via monitoring functions. \emph{52nd IEEE Conference on
  Decision and Control}, pages 1031--1036, 2013.

\bibitem{A:1957} T. Apostol, {\it Mathematical Analysis – A Modern Approach to Advanced Calculus}, Massachusetts: Addison-Wesley Publishing Company, 1957.

\bibitem{s8}  
K.~J.~{\r A}str{\" o}m and B.~P.~Bernhardsson, 
\newblock Comparison of periodic and event based sampling for first-order stochastic systems.
\newblock {\em IFAC World Congress}, 32:5006--5011, 1999.


\bibitem{paper20}
A.~Battistel, T.~R. Oliveira, V.~H.~P. Rodrigues, and L.~Fridman,
  Multivariable binary adaptive control using higher-order sliding modes
  applied to inertially stabilized platforms. \emph{European Journal of
  Control}, 63:28--39, 2022.


\bibitem{BH:2013}
D.~P.~Borgers and W.~P. M.~H.~Heemels, 
\newblock On minimum inter-event times in event-triggered control.
\newblock In {\em IEEE 52nd  IEEE Conf. Decis. Control.}, pages 7370--7375, Firenze, Italy, 2013.

\bibitem{CKAL:2002}
J.-Y.~Choi, M.~Krsti{\'c}, K.~B.~Ariyur, and J.~S.~Lee, 
\newblock Extremum seeking control for discrete-time systems.
\newblock {\em IEEE Trans. Automat. Contr.}, 47:318--323, 2002.

\bibitem{source1}
J. Cochran and M. Krstic, 
\newblock Nonholonomic source seeking with tuning of angular velocity.
\newblock {\em IEEE Transactions on Automatic Control}, 54:717--731, 2009.

\bibitem{paper5}
C.~L. Coutinho, T.~R. Oliveira, and J.~P. V.~S. Cunha, Output-feedback
  sliding-mode control via cascade observers for global stabilisation of a
  class of nonlinear systems with output time delay. \emph{International
  Journal of Control}, 87(11):2327--2337, 2014.

\bibitem{FKB:2012}
P.~Frihauf, M.~Krsti{\' c}, and T.~Ba{\c s}ar, 
\newblock Nash equilibrium seeking in non-cooperative games.
\newblock {\em IEEE Trans. Automat. Contr.}, 57:1192--1207, 2012.

\bibitem{GKN:2012}
A.~Ghaffari, M.~Krsti{\'c}, and D.~Ne{\u s}ic, 
\newblock Multivariable {N}ewton-based extremum seeking.
\newblock {\em Automatica}, 48:1759--1767, 2012.

\bibitem{G:2014} A. Girard, "Dynamic triggering mechanism for event-triggered control," {\it IEEE Trans. Autom. Control}, vol. 60, pp. 1992--1997, 2014.

\bibitem{paper4}
L.~L. Gomes, L.~Leal, T.~R. Oliveira, J.~P. V.~S. Cunha, and T.~C. Revoredo,
  ``Unmanned quadcopter control using a motion capture system,'' \emph{IEEE
  Latin America Transactions}, vol.~14, no.~8, pp. 3606--3613, 2016.

\bibitem{Basar:2019}
Z. Han, D. Niyato, W. Saad, and T. Başar, Game Theory for Next Generation Wireless and Communication Networks: Modeling, Analysis, and Design, Cambridge University Press, 2019.

\bibitem{s5}  
W.~P.~M.~H.~Heemels, M.~C.~F.~Donkers, and A.~R.~Teel, 
\newblock Periodic event-triggered control for linear systems.
\newblock {\em IEEE Trans. Automat. Contr.}, 58:847--861, 2012.

\bibitem{HJT:2012} W. P. M. H. Heemels, K. H. Johansson, and P. Tabuada, "An introduction to event-triggered and self-triggered control," {\it IEEE Conf. Decis. Control}, pp. 3270--3285, 2012.

\bibitem{HNX:2007}
J.~P.~Hespanha, P.~Naghshtabrizi, and Y.~Xu, 
\newblock A survey of recent results in networked control systems.
\newblock {\em Proceedings of the IEEE}, 95(1):138--162.

\bibitem{Johansson} 
W. Huo, K. F. E. Tsang, Y. Yan, K. H. Johansson, and L. Shi,  
\newblock Distributed Nash equilibrium seeking with stochastic event-triggered mechanism. 
\newblock {\em Automatica}, vol. 162, paper no. 111486, 2024.

\bibitem{K:2002} H. K. Khalil, {\it Nonlinear Systems}, Upper Saddle River, NJ: Prentice Hall, 2002.

\bibitem{KW:2000}
M.~Krsti{\' c} and H.-H.~Wang, 
\newblock Stability of extremum seeking feedback for general nonlinear dynamic systems.
\newblock {\em Automatica}, 36:595--601, 2000.

\bibitem{LK:2010}
S.-J.~Liu and M.~Krsti{\'c}, 
\newblock Stochastic averaging in continuous time and its applications to extremum seeking.
\newblock {\em IEEE Trans. Automat. Contr.}, 55:2235--2250, 2010.

\bibitem{MK:2009}
C.~Manzie and M.~Krsti{\'c}, 
\newblock Extremum seeking with stochastic perturbations.
\newblock {\em IEEE Trans. Automat. Contr.}, 54:580--585, 2009.

\bibitem{Mazo_Tabuada} 
M. Mazo and P. Tabuada, Decentralized event-triggered control over wireless sensor/actuator networks. IEEE Transactions on Automatic Control 56(10):2456--2461, 2011.

\bibitem{s9b}
S.~Monaco and D.~Normand-Cyrot, 
\newblock Discrete time models for robot arm control.
\newblock {\em IFAC Proceedings Volumes}, 18:525--529, 1985.


\bibitem{paper14}
T.~R. Oliveira, L.~R. Costa, J.~M.~Y. Catunda, A.~V. Pino, W.~Barbosa, and
  M.~N. de~Souza, Time-scaling based sliding mode control for neuromuscular
  electrical stimulation under uncertain relative degrees. \emph{Medical
  Engineering $\&$ Physics}, 44:53--62, 2017.

\bibitem{paper12}
T.~R. Oliveira, J.~P. V.~S. Cunha, and L.~Hsu, Adaptive sliding mode control
  based on the extended equivalent control concept for disturbances with
  unknown bounds. \emph{Advances in Variable Structure Systems and Sliding
  Mode Control---Theory and Applications. Studies in Systems, Decision and
  Control}, 115:149--163, 2017.


\bibitem{paper17}
T.~R. Oliveira, L.~Hsu, and A.~J. Peixoto, Output-feedback global tracking
  for unknown control direction plants with application to extremum-seeking
  control. \emph{Automatica}, 47(9):2029--2038, 2011.


\bibitem{paper9}
T.~R. Oliveira and M.~Krstic, Newton-based extremum seeking under actuator
  and sensor delays. \emph{IFAC-PapersOnLine}, 48(12):304--309,
  2015.

\bibitem{OK:2022_book_golden}
T.~R.~Oliveira and M.~Krstic, 
\newblock {\em Extremum Seeking through Delays and PDEs}.
\newblock SIAM, Philadelphia, 2022.

\bibitem{OKT:2017}
T.~R.~Oliveira, M.~Krsti{\' c}, and D.~Tsubakino, 
\newblock Extremum seeking for static maps with delays.
\newblock {\em IEEE Trans. Automat. Contr.}, 62:1911--1926, 2017.


\bibitem{paper10}
T.~R. Oliveira, A.~C. Leite, A.~J. Peixoto, and L.~Hsu, Overcoming
  limitations of uncalibrated robotics visual servoing by means of sliding mode
  control and switching monitoring scheme. \emph{Asian Journal of Control},
  16(3):752--764, 2014.


\bibitem{paper11}
T.~R. Oliveira, A.~J. Peixoto, and L.~Hsu, Peaking free output‐feedback
  exact tracking of uncertain nonlinear systems via dwell‐time and norm
  observers. \emph{International Journal of Robust and Nonlinear Control},
  23(5):483--513, 2013.
	
	\bibitem{paper15}
T.~R. Oliveira, A.~J. Peixoto, and L.~Hsu, Global tracking for a class of
  uncertain nonlinear systems with unknown sign-switching control direction by
  output feedback. \emph{International Journal of Control}, 88(9):1895--1910, 2015.


\bibitem{paper3}
T.~R. Oliveira, A.~J. Peixoto, and E.~V.~L. Nunes, Binary robust adaptive
  control with monitoring functions for systems under unknown
  high‐frequency‐gain sign, parametric uncertainties and unmodeled
  dynamics. \emph{International Journal of Adaptive Control and Signal
  Processing}, 30(8-10):1184--1202, 2016.
	
	
	\bibitem{paper18}
T.~R. Oliveira, A.~J. Peixoto, E.~V.~L. Nunes, and L.~Hsu, Control of
  uncertain nonlinear systems with arbitrary relative degree and unknown
  control direction using sliding modes. \emph{International Journal of
  Adaptive Control and Signal Processing}, 21(8-9):692--707, 2007.
	
		
\bibitem{paper1}
T.~R. Oliveira, V.~H.~P. Rodrigues, and L.~Fridman, Generalized model
  reference adaptive control by means of global {HOSM} differentiators. 
  \emph{IEEE Transactions on Automatic Control}, 64(5):2053--2060, 2018.


\bibitem{paper16}
T.~R. Oliveira, V.~H.~P. Rodrigues, M.~Krstic, and T.~Basar, Nash equilibrium
  seeking in quadratic noncooperative games under two delayed
  information-sharing schemes. \emph{Journal of Optimization Theory and
  Applications}, 191(2):700--735, 2021.


\bibitem{paper7}
A.~J. Peixoto, T.~R. Oliveira, L.~Hsu, F.~Lizarralde, and R.~R. Costa, Global
  tracking sliding mode control for a class of nonlinear systems via variable
  gain observer. \emph{International Journal of Robust and Nonlinear
  Control}, 21(2):177--196, 2011.


\bibitem{paper13}
H.~L. C.~P. Pinto, T.~R. Oliveira, and L.~Hsu, Sliding mode observer for
  fault reconstruction of time-delay and sampled-output systems---a time shift
  approach. \emph{Automatica}, 106:390--400, 2019.


\bibitem{P:1979} V. A. Plotnikov, "Averaging of differential inclusions," {\it Ukrainian Mathematical Journal}, vol. 31, pp. 454--457, 1980.

\bibitem{paper8}
V.~H.~P. Rodrigues and T.~R. Oliveira, Global adaptive {HOSM} differentiators
  via monitoring functions and hybrid state-norm observers for output
  feedback. \emph{International Journal of Control}, 91(9):2060--2072, 2018.

\bibitem{AUT:2025} V. H. P. Rodrigues, T. R. Oliveira, L. Hsu, M. Diagne, and M. Krsti{\' c}, "Event-triggered and periodic event-triggered extremum seeking control," {\it Automatica}, vol. 174, paper no. 112161, pp. 1--16, 2025.

\bibitem{paper2}
D.~Rusiti, G.~Evangelisti, T.~R. Oliveira, M.~Gerdts, and M.~Krstic,
  Stochastic extremum seeking for dynamic maps with delays. \emph{IEEE
  Control Systems Letters}, 3(1):61--66, 2019.

\bibitem{source2}
R. Suttner and M. Krstic,
\newblock Source seeking with a torque-controlled unicycle.
\newblock {\em  IEEE Control Systems Letters}, 7:79--84, 2023.

\bibitem{source4}
R. Suttner and M. Krstic,
\newblock Nonlocal nonholonomic source seeking despite local extrema.
\newblock {\em   IEEE Transactions on Automatic Control}, 69:2575--2582, 2024.

\bibitem{source5}
R. Suttner and M. Krstic, 
\newblock Overcoming local extrema in torque-actuated source seeking using the divergence theorem and delay.
\newblock {\em Automatica}, paper no. 111799, 2024. 


\bibitem{T:2007}
P.~Tabuada, 
\newblock Event-triggered real-time scheduling of stabilizing control tasks.
\newblock {\em IEEE Trans. Automat. Contr.}, 52:1680--1685, 2007.

\bibitem{source6}
V. Todorovski and M. Krstic,
\newblock Newton nonholonomic source seeking for distance-dependent maps.
\newblock {\em IEEE Trans. Automat. Contr.}, 70:510--517, 2025.

\bibitem{paper19}
D.~Tsubakino, T.~R. Oliveira, and M.~Krstic, Extremum seeking for distributed
  delays. \emph{Automatica}, 153(111044):1--14, 2023.

\bibitem{Wang_Lemmon} 
X. Wang and M. D. Lemmon, Event-triggering in distributed networked control systems. IEEE Transactions on Automatic Control, 56(3):586--601, 2011.

\bibitem{source3}
B. Wang, S. Nersesov, H. Ashrafiouon, P. Naseradinmousavi, and M. Krstic, \newblock Underactuated source seeking by surge force tuning: Theory and boat experiments.
\newblock {\em IEEE Transactions on Control Systems Technology}, 31: 1649--1662, 2023.

\bibitem{W:1971} F. Warner, {\it Foundations of Differentiable Manifolds and Lie Groups}, Chicago, IL: Scott Foresman and Company, 1971.


\bibitem{s1} 
J.~K.~Yook, D.~M.~Tilbury, and N.~R.~Soparkar,    
\newblock Trading computation for bandwidth: Reducing communication in distributed control systems using state estimators.
\newblock {\em IEEE transactions on Control Systems Technology}, 10:503--518, 2002.

\bibitem{SCL_2007}
C. Zhang, D. Arnold, N. Ghods, A. Siranosian, and M. Krstic,
\newblock Source seeking with nonholonomic unicycle without position measurement and with tuning of forward velocity. 
\newblock {\em Systems and Control Letters}, 56:245--252, 2007.

\bibitem{ZHGDDYP:2020}
X.-M.~Zhang, Q.-L.~Han, X.~Ge, D.~Ding, L.~Ding, D.~Yue, and C.~Peng, 
\newblock Networked control systems: A survey of trends and techniques.
\newblock {\em IEEE/CAA J. Autom. Sin.}, 7(1):1--17, 2020.

\bibitem{ZFO:2023}
Y. Zhu, E. Fridman, and T. R. Oliveira, 
\newblock Sampled-data extremum seeking with constant delay: a time-delay approach.
\newblock {\em IEEE Trans. Automat. Contr.}, 68:432--439, 2023. 

\end{thebibliography}
\end{document}